\documentclass[journal]{IEEEtran}

\usepackage{amsmath}
\usepackage{amsfonts}
\usepackage{amssymb}
\usepackage{amsthm}

\usepackage{mathdots} 

\usepackage{array}

\usepackage{graphicx}
\usepackage{hyperref} 

\usepackage{color}
\usepackage{import}

\usepackage{parskip}

\newcommand{\reals}{\mathbb{R}}
\newcommand{\defeq}{\mathrel{\mathop:}=}

\newcommand{\A}{\mathcal{A}}
\newcommand{\dual}{\widetilde}

\newcommand\scalemath[2]{\scalebox{#1}{\mbox{\ensuremath{\displaystyle #2}}}}

\IEEEoverridecommandlockouts
\usepackage{tikz}
\usepackage{textcomp}
\usepackage{hyperref}
\usepackage{lipsum}

\newcommand\copyrighttext{%
  \footnotesize \textcopyright 2016 IEEE. Personal use of this material is permitted.
  Permission from IEEE must be obtained for all other uses, in any current or future
  media, including reprinting/republishing this material for advertising or promotional
  purposes, creating new collective works, for resale or redistribution to servers or
  lists, or reuse of any copyrighted component of this work in other works.
  DOI: \url{https://doi.org/10.1109/MSP.2016.2594277}}
\newcommand\copyrightnotice{%
\begin{tikzpicture}[remember picture,overlay]
\node[anchor=south,yshift=3pt] at (current page.south) {\fbox{\parbox{\dimexpr\textwidth-\fboxsep-\fboxrule\relax}{\copyrighttext}}};
\end{tikzpicture}%
}

\begin{document}
\title{Efficient Adjoint Computation for Wavelet and Convolution Operators}
\author{James~Folberth and 
        Stephen~Becker%
\thanks{J. Folberth and S. Becker are with the Department
of Applied Mathematics, University of Colorado at Boulder,
Boulder, CO, 80309, USA.}}

\markboth{IEEE Signal Processing Magazine: Lecture Notes}%
{J. Folberth and S. Becker: Fast Adjoint Wavelet}

\maketitle
\copyrightnotice

\begin{abstract}
   First-order optimization algorithms, often preferred for large problems, require the gradient of the differentiable terms in the objective function.  These gradients often involve linear operators and their adjoints, which must be applied rapidly.  We consider two example problems and derive methods for quickly evaluating the required adjoint operator.  The first example is an image deblurring problem, where we must compute efficiently the adjoint of multi-stage wavelet reconstruction.  Our formulation of the adjoint works for a variety of boundary conditions, which allows the formulation to generalize to a larger class of problems.  The second example is a blind channel estimation problem taken from the optimization literature where we must compute the adjoint of the convolution of two signals.  In each example, we show how the adjoint operator can be applied efficiently while leveraging existing software.
\end{abstract}

\subsection*{Prerequisites}
The reader should be familiar with linear algebra, wavelets, and basic Fourier analysis.  Knowledge of first-order iterative optimization algorithms is beneficial for motivating the need for a fast adjoint computation, but should not be necessary to understand the adjoint computation itself.

\subsection*{Motivation}
Many estimation problems are modeled using the composition of a differentiable loss $\ell$ with an affine map $x\mapsto\mathcal{A}x+b$, the simplest such examples being generalized linear models such as least-square data fitting. The gradient of $\ell \circ \A$ is given by the chain rule as $\mathcal{A}^* \nabla \ell( \A x + b)$. For small problems, $\A$ may be explicitly represented as a matrix and $\A^*$ calculated as the conjugate transpose of the matrix. 

In contrast, larger problems are often carefully formulated to involve only linear operators $\A$ which have a fast transform (e.g., FFT, wavelet transform, convolution) and are therefore never explicitly stored as matrices.  While many software packages provide fast routines to compute $\A^{-1}$ or the pseudoinverse $\A^\dagger$, few packages include routines for $\A^*$.

The purpose of this note is to describe a few cases where one can compute the action of $\A^*$ with the same complexity as applying $\A$. Boundary conditions play a confounding role in the calculation, and we present a simple framework to correctly take them into account.

\section{Example 1: An Image Deblurring Problem}
Digital images can be blurred through a variety of means.  For example, the optical system can be out of focus and/or atmospheric turbulence can cause blurring of astronomical images.  The goal of image deblurring is to recover the original, sharp image. One class of techniques poses the deblurring problem  as an optimization problem, where the optimization variables are the wavelet coefficients of the recovered image.  We assume  that the blurring operator $\mathcal{R}$ is known. For instance, $\mathcal{R}$ could represent a Gaussian point spread function (PSF) under symmetric (reflexive) boundary conditions. The blurring operator is usually singular or ill-conditioned \cite{hansen_2006}.

Let $\mathcal{W}$ represent a multi-level wavelet reconstruction operator with suitable boundary conditions, $b$ be the observed, blurry image, and $\mathcal{A}=\mathcal{RW}$ be the linear operator that synthesizes an image from wavelet coefficients $x$ and then blurs the synthesized image under the blurring operator $\mathcal{R}$.  In order to both regularize the problem and to take into account that many real-world images have a sparse wavelet representation, we solve the following standard model which seeks to recover sparse wavelet coefficients that accurately reconstruct the blurred image:

\begin{equation}
\label{eq:syn_problem}
\min_x~ \dfrac{1}{2}\|\mathcal{RW}x-b\|_2^2 + \lambda \|x\|_1, \quad \lambda>0.
\end{equation}
This commonly used formulation \cite{beck_2009} is sometimes called the synthesis formulation, since we reconstruct an image $\mathcal{W}x$ from the wavelet coefficients in the optimization variable $x$.  We seek to find an $x$ that accurately reconstructs the blurred image $b$ while simultaneously having only a few nonzero entires.

To clarify our notation, suppose we have an image $b$.  In terms of the wavelet reconstruction operator $\mathcal{W}$, we compute the wavelet coefficients of $b$ with $x=\mathcal{W}^\dagger b$, which corresponds to wavelet analysis.  In each case, $\mathcal{W}$ and $\mathcal{W}^\dagger$ implicitly impose and handle the boundary conditions we impose on $b$ (e.g., zero, periodic, symmetric).

The gradient of the differentiable term $f(x)={\frac{1}{2}\|\mathcal{RW}x-b\|_2^2}$ is

\[ \nabla f(x) = \mathcal{W}^\ast \mathcal{R}^\ast(\mathcal{RW}x-b). \] 

\noindent In the case of a blurring PSF, we can apply $\mathcal{R}$ and $\mathcal{R}^\ast$ in $\mathcal{O}(N\log N)$ time in the Fourier domain via the FFT \cite{beck_2009, hansen_2006}.  Wavelet software packages (e.g., {\sc matlab} Wavelet Toolbox \cite{matlab_wt_2015}) implement discrete wavelet analysis and reconstruction in $\mathcal{O}(N)$ time; these operations correspond to $\mathcal{W}^\dagger$ and $\mathcal{W}$.  However, the adjoint of wavelet reconstruction, $\mathcal{W}^\ast$, is not implemented in software packages and is not expressly mentioned in the literature \cite{mallat_2009, strang_1996}.  A ``trick'' commonly used in practice is the ``pseudoinverse approximation'' $\mathcal{W}^\dagger\approx\mathcal{W}^\ast$, but we would like to find a way of computing the true adjoint of wavelet reconstruction (and analysis) in the same $\mathcal{O}(N)$ time.

To derive the true adjoint, we first separate $\mathcal{W}$ into signal extension and wavelet transform.  We recall a few common signal extensions, consider the case of orthogonal wavelets, and then use a property of frames to construct the adjoint for biorthogonal wavelets.  It transpires that we only need to implement the adjoint of the signal extensions and may use existing software for the wavelet transforms.

\subsection{A few signal extensions}
To implement a wavelet transform with boundary conditions, a popular method is to extend the given signal to satisfy the desired boundary conditions and then transform the extended signal; this is the approach used in {\sc matlab}'s Wavelet Toolbox \cite{matlab_wt_2015}.  If we let $\mathcal{E}$ be the desired extension operator and $\mathcal{W}^\dagger_\text{zpd}$ be the wavelet analysis operator for extended signals under zero boundary conditions, we can write wavelet analysis as 

\[ \mathcal{W}^\dagger = \mathcal{W}^\dagger_\text{zpd}\mathcal{E}. \] 

Consider a $1$-dimensional signal $y[n]$, $n=0,...,N-1$.  \footnote{We consider only $1$d signals here for brevity.  The results can be extended to 2d quite naturally.}  Let $L$ be the maximum length of the wavelet analysis filters.  The double sided zero-padded extension of $y[n]$ puts $L'\defeq L-1$ zeros on each end of the signal:
\[ \underbrace{0, ..., 0}_{L'}, y[0], ..., y[N-1], \underbrace{0, ..., 0}_{L'}. \]

\noindent We can write down the linear operator $\mathcal{E}_\text{zpd}$ as a matrix that performs this operation on $y[n]$.

\[ \mathcal{E}_\text{zpd} = \begin{bmatrix} 0_{L'\times N}\\ I_{N\times N}\\ 0_{L'\times N}\end{bmatrix}. \] 

   \noindent We see that $\mathcal{E}_\text{zpd}$ has orthonormal columns, so that ${\mathcal{E}_\text{zpd}^\dagger\mathcal{E}_\text{zpd} = \mathcal{E}_\text{zpd}^\ast\mathcal{E}_\text{zpd}=I}$.

   Note that $\mathcal{W}_\text{zpd}$ implemented in existing software may also implicitly contain a zero-padded extension.  When we apply wavelet analysis, $\mathcal{W}_\text{zpd}^\dagger$, we implicitly apply $\mathcal{E}_\text{zpd}^\dagger$.  However, since zero-padding is orthonormal, $\mathcal{E}_\text{zpd}^\dagger=\mathcal{E}_\text{zpd}^\ast$, and so we can, for a general wavelet transform, write the adjoint as 

\begin{equation}
\label{eq:wzpd_adjoint}
\mathcal{W}^\ast = \mathcal{W}_\text{zpd}^\ast\left(\mathcal{E}^\dagger\right)^\ast.
\end{equation}

\noindent The splitting in equation (\ref{eq:wzpd_adjoint}) allows us to correctly handle the adjoint extension ourselves.  It remains to implement $\left(\mathcal{E}^\dagger\right)^\ast$ and determine how to apply $\mathcal{W}^\ast_\text{zpd}$ with existing software.

The zero-padded extension is not often used by itself since it can introduce boundary artifacts.
 A better extension mode is the half-point symmetric extension, which is the default extension mode in {\sc matlab}'s Wavelet Toolbox \cite{strang_1996, matlab_wt_2015}.  The double sided half-point symmetric extension reflects the signal about its boundaries in the following manner:

\[ \underbrace{y[L-1], ..., y[0]}_\text{Left extension}, y[0], ..., y[N-1], \underbrace{y[N-1], ..., y[N+L-2]}_\text{Right extension}. \] 

\noindent As in the zero-padded case, we can readily form a matrix that performs this operation on $y[n]$.  Let $P$ be the appropriately sized permutation matrix that reverses the ordering of columns.  Then we have

\[ \scalemath{0.8}{\mathcal{E}_\text{sym} = \begin{bmatrix} & & \iddots & & &\\ & 1 &&&&\\ 1&&&&&\\1&&&&\\&1&&&&\\&&\ddots&&&\\&&&1&\\&&&&1\\&&&&1\\&&&1&\\&&\iddots&&\end{bmatrix} = \begin{bmatrix}I_{L'\times L'}P & 0 & 0\\I_{L'\times L'} & 0 & 0\\0 & I_{N-2L'\times N-2L'} & 0\\0&0&I_{L'\times L'}\\0&0&I_{L'\times L'}P\end{bmatrix}.} \]

\noindent The adjoint of the pseudoinverse can be computed in closed form and essentially amounts to rescaling the nonzero entries of $\mathcal{E}_\text{sym}$.  Assuming $N > 2L'$, which usually occurs in practice, the adjoint of the  pseudoinverse of $\mathcal{E}_\text{sym}$ is 

\[ \left(\mathcal{E}_\text{sym}^\dagger\right)^\ast = \begin{bmatrix}\frac{1}{2}I_{L'\times L'}P & 0 & 0\\\frac{1}{2}I_{L'\times L'} & 0 & 0\\0 & I_{N-2L'\times N-2L'} & 0\\0&0&\frac{1}{2}I_{L'\times L'}\\0&0&\frac{1}{2}I_{L'\times L'}P\end{bmatrix}. \]

\noindent If ${N \le 2L'}$, the form of the pseudoinverse is slightly different, with some of the $1/2$ terms becoming $1/3$.

Another standard extension is the periodic extension.  Similar to the half-point symmetric extension, the periodic extension operator is

\[ \mathcal{E}_\text{per} = \begin{bmatrix} 0 & 0 & I_{L'\times L'}\\&I_{N\times N}&\\I_{L'\times L'}&0&0\end{bmatrix} \] 

\noindent and the adjoint of its pseudoinverse is

\[ \left(\mathcal{E}_\text{per}^\dagger\right)^\ast = \begin{bmatrix} 0 & 0 & \frac{1}{2}I_{L'\times L'}\\\frac{1}{2}I_{L'\times L'}&0&0\\0&I_{N-2L'\times N-2L'}&0\\0&0&\frac{1}{2}I_{L'\times L'}\\\frac{1}{2}I_{L'\times L'}&0&0\end{bmatrix}. \] 

\noindent Again, in the unusual case that $N\le 2L'$, the form of the pseudoinverse changes slightly.  By viewing other signal extension operators as matrices acting on signal vectors, one can readily find the appropriate $(\mathcal{E}^\dagger)^\ast$.

\subsection{The adjoint for orthogonal wavelets}
For orthogonal wavelets (e.g., Haar and more generally Daubechies wavelets), the adjoint $\mathcal{W}_\text{zpd}^\ast$ is exactly the analysis operator $\mathcal{W}_\text{zpd}^\dagger$.  For orthogonal wavelets with general boundary conditions, we can use the splitting (\ref{eq:wzpd_adjoint}) and implement $\left(\mathcal{E}^\dagger\right)^\ast$.

\begin{figure}
   \centering
   \includegraphics[width=0.8\columnwidth]{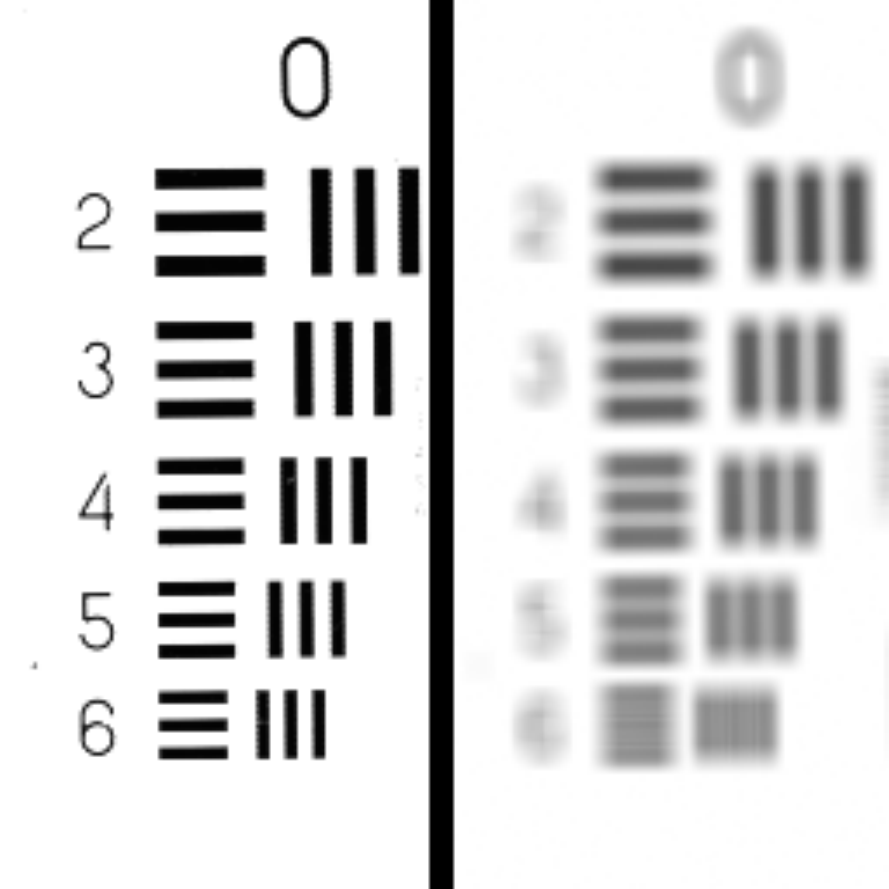}
   \caption{The left half of the original, unblurred image and the left half of the blurred image.  To create the blurred image $b$, the blurring operator $\mathcal{R}$ is applied to the original image and small amount of Gaussian noise is added.}
   \label{fig:original}
\end{figure}

For an example image, consider a standard resolution test chart \cite{weber_1993}, shown in Figure \ref{fig:original}.  We prescribe symmetric boundary conditions on the image, which will affect both $\mathcal{W}$, $\mathcal{R}$, and their adjoints.  Symmetric boundary conditions (compared to periodic or zero-padded boundary conditions) produce significantly fewer edge effects, and so are a natural choice.

Fortuitously, for orthogonal wavelets, the adjoint of $\mathcal{W}$ with symmetric boundary conditions is ``very close to'' the analysis operator $\mathcal{W}^\dagger$.  This fact is often used (e.g., \cite{beck_2009}), and as we can see from the splitting (\ref{eq:wzpd_adjoint}), the error is introduced through $\mathcal{E} \approx \left(\mathcal{E}^\dagger\right)^\ast$.  Again, by merely implementing $\left(\mathcal{E}^\dagger\right)^\ast$, we can use the true adjoint $\mathcal{W}^\ast$.

\begin{figure}
   \centering
   \includegraphics[width=0.8\columnwidth]{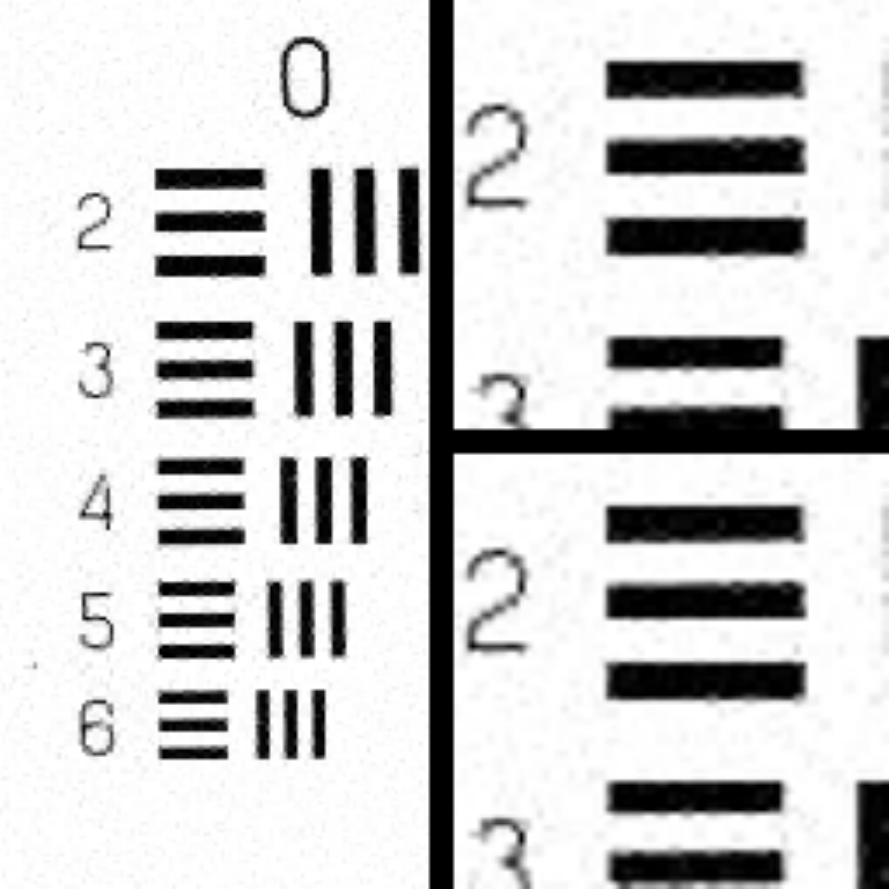}
   \caption{Deblurred image after 2500 iterations of FISTA with $\mathcal{W}$ a three-stage Haar transform with symmetric boundary conditions.  The left half of the figure shows the deblurred image using the true adjoint $\mathcal{W}^\ast$.  The top-right shows a zoomed-in view using $\mathcal{W}^\ast$; the bottom-right shows the same zoomed-in view using $\mathcal{W}^\dagger$ and 2500 iterations of FISTA.}
   \label{fig:compare_db1_sym}
\end{figure}

Figure \ref{fig:compare_db1_sym} shows the deblurred image after 2500 iterations of FISTA, a fast proximal gradient method, developed by Beck and Teboulle \cite{beck_2009}.  The wavelet reconstruction operator $\mathcal{W}$ is taken to be a three-stage Haar discrete wavelet transform with symmetric boundary conditions.  We set $\lambda=2\times 10^{-5}$.  We use both the true adjoint $\mathcal{W}^\ast$ and the approximation $\mathcal{W}^\dagger$.  Using the true adjoint, the image reconstruction relative error (vs. the unblurred image) is $8.91\times10^{-4}$ and $31.8\%$ of the coefficients are nonzero.  Using the pseudoinverse approximation, the relative error is $8.91\times 10^{-4}$ and $32.1\%$ of the coefficients are nonzero.  In this case the pseudoinverse approximation works extremely well.

\subsection{The adjoint for biorthogonal wavelets}
For biorthogonal wavelets, we no longer have $\mathcal{W}_\text{zpd}^\ast=\mathcal{W}_\text{zpd}^\dagger$.  Since we have relaxed ourselves to biorthogonal wavelets (which are nice for symmetric boundary conditions \cite{rout_2003}), it is perhaps too much to ask that the adjoint of the primal wavelet reconstruction operator involves only the primal wavelets.  Let us recall briefly the pertinent aspects of frames of $\reals^d$.  These facts are described in more detail and generality in \cite{mallat_2009}.

Let $y\in\reals^N$ be an arbitrary signal vector and $\{\phi_i\}_{i=1}^p$, $p\ge N$, be a set of vectors in $\reals^N$.  Define the analysis operator $\Phi$ as the $p\times N$ matrix

\[ \Phi = \begin{bmatrix}\phi_1^T\\\vdots\\\phi_p^T\end{bmatrix}. \] 

   \noindent If $\Phi$ has full rank, we say that $\{\phi_i\}$ is a frame and we call $\Phi$ the frame analysis operator.  Henceforth we assume $\{\phi_i\}$ is a frame.
   
   The product $\Phi u$ computes the expansion coefficients of the signal ${u\in\reals^N}$ in the frame $\{\phi_i\}$.  The frame synthesis operator $\Phi^\ast$ constructs a vector in $\reals^N$ given some expansion coefficients.  Since $\{\phi_i\}$ is assumed to be a frame, $\Phi^\ast\Phi$ is invertible and we may define the Moore-Penrose pseudoinverse $\Phi^\dagger$, which implements reconstruction in the frame as $\Phi^\dagger=\left(\Phi^\ast\Phi\right)^{-1}\Phi$.

   A dual frame can be associated with a (primal) frame by defining the dual frame vectors $\dual{\phi}_i = \left(\Phi^\ast\Phi\right)^{-1}\phi_i$ for $i=1,...,p$.  We may define analogously the dual frame analysis operator $\dual{\Phi}$ and dual frame synthesis operator $\dual{\Phi}^\ast$.

   \noindent A fundamental relationship between primal and dual frame operators is $\dual{\Phi}^\ast = \Phi^\dagger$ (see \cite{mallat_2009}, Theorem 5.5).  We can directly relate the frame operators and the wavelet operators under zero-padded boundary conditions.  We can write zero-padded wavelet reconstruction as frame reconstruction: $\mathcal{W}_\text{zpd}=\Phi^\dagger$.  Then using the above frame relations, we have

   \begin{equation}
      \label{eq:wave_dual_frame}
      \mathcal{W}_\text{zpd}^\ast = \left(\Phi^\dagger\right)^\ast = \left(\Phi^\ast\right)^\dagger = \dual{\Phi}.
   \end{equation}

\begin{table}[ht]
   \centering
   \caption{A few primal and dual frame relationships.}
   \label{tab:frame_relations}
   \begin{tabular}{lll}
      &Primal frame&Dual frame\\\hline\\[-0.5em]
      Analysis & $\Phi = \mathcal{W}_\text{zpd}^\dagger$ & $\dual{\Phi}=\left(\Phi^\dagger\right)^\ast$\\\\[-0.5em]
      Synthesis & $\Phi^\ast = \left(\mathcal{W}_\text{zpd}^\dagger\right)^\ast$ & $\dual{\Phi}^\ast = \Phi^\dagger$ \\\\[-0.5em]
      Reconstruction & $\Phi^\dagger = \mathcal{W}_\text{zpd}$ & $\dual{\Phi}^\dagger = \Phi^\ast$
   \end{tabular}
\end{table}

   \noindent We present other relationships between the primal and dual frames in Tabel \ref{tab:frame_relations}.  Note that $\Phi^\dagger\Phi=I$ but $\Phi\Phi^\dagger\neq I$ in general.

   In (\ref{eq:wave_dual_frame}) we have $\mathcal{W}_\text{zpd}^\ast$ in terms of dual frame analysis.  For biorthogonal wavelets, dual frame analysis corresponds to wavelet analysis with the dual wavelets, which has a fast transform and is typically implemented together with the primal wavelets \cite{matlab_wt_2015}.  For biorthogonal (and orthogonal) wavelets, this states $\mathcal{W}^\ast_\text{zpd}=\dual{\mathcal{W}}_\text{zpd}^\dagger$, the latter of which is a standard operation and has a fast implementation.  Combining this with the signal extension operator, we finally find

   \begin{equation}
   \label{eq:full_wave_adjoint}
      \mathcal{W}^\ast = \mathcal{W}_\text{zpd}^\ast\left(\mathcal{E}^\dagger\right)^\ast = \dual{\mathcal{W}}_\text{zpd}^\dagger\left(\mathcal{E}^\dagger\right)^\ast.
   \end{equation}
   
   \noindent From this we see that with a complete wavelet software library, one only needs to implement $(\mathcal{E}^\dagger)^\ast$ to use the true adjoint instead of the approximation $\mathcal{W}^\ast\approx\mathcal{W}^\dagger$.



Consider deblurring the resolution chart from Figure \ref{fig:original} with 2500 iterations of FISTA using a three-stage CDF 9/7 discrete wavelet transform with symmetric boundary conditions.  The results of deblurring are visually very similar to Figure \ref{fig:compare_db1_sym}.  Using the true adjoint, the image reconstruction relative error (vs. the unblurred image) is $9.45\times 10^{-4}$ and $29.23\%$ of the coefficients are nonzero.  Using the pseudoinverse approximation, the relative error is $9.67\times 10^{-4}$ and $29.74\%$ of the coefficients are nonzero.  Here the pseudoinverse approximation works, but not quite as well as the true adjoint.

\begin{figure}
   \centering
   \includegraphics[width=\columnwidth]{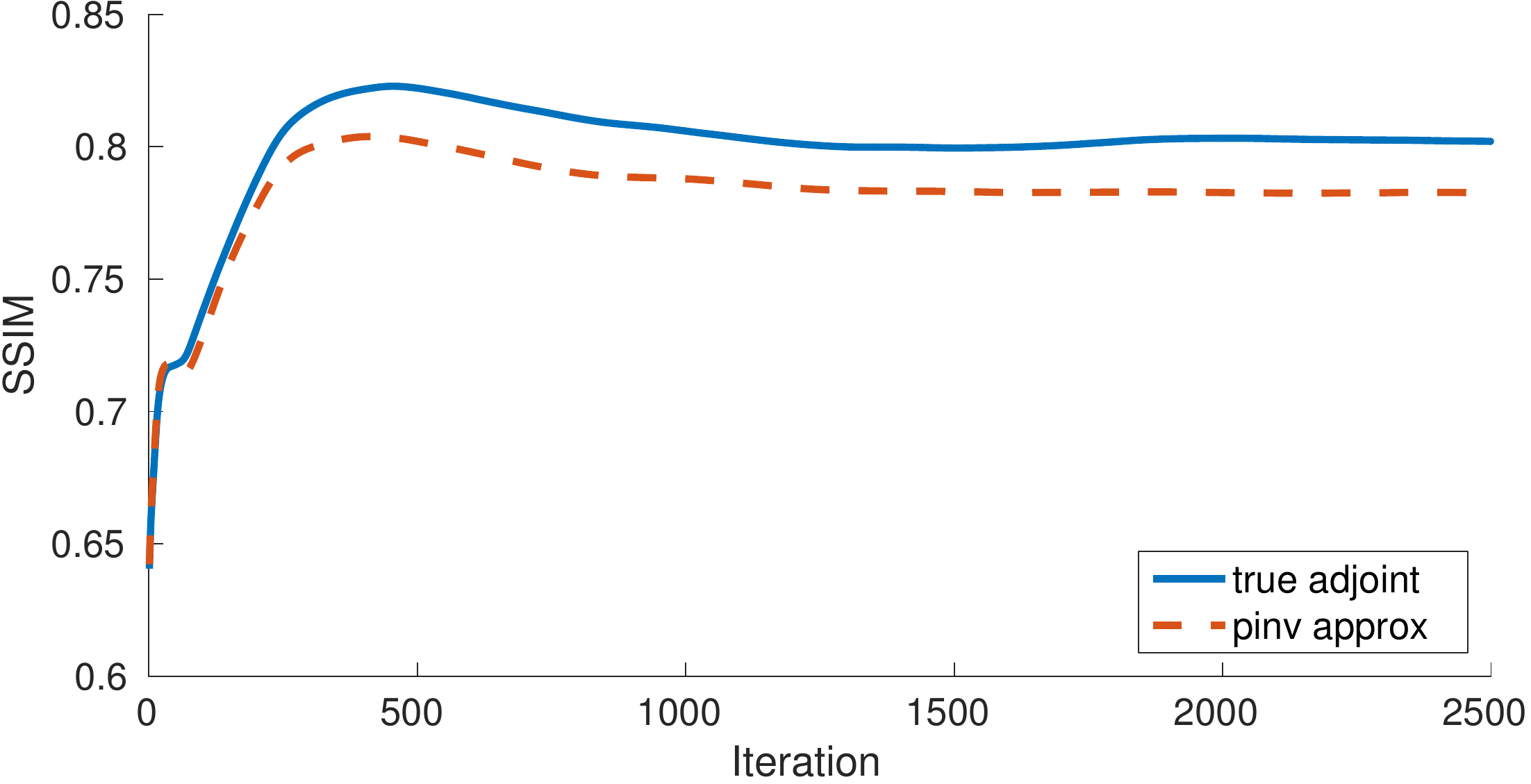}
   \caption{Structural similarity (SSIM) index over 2500 iterations of FISTA with $\mathcal{W}$ a three-stage CDF 9/7 transform and using the true adjoint and the pseudoinverse approximation.  SSIM compares the deblurred image to the original, unblurred image.}
   \label{fig:adjoint_convergence}
\end{figure}

Figure \ref{fig:adjoint_convergence} shows the structural similarity index \cite{wang_2004} for the CDF 9/7 experiment using both the pseudoinverse approximation and true adjoint.  For the first few hundred iterations of FISTA, the pseudoinverse approximation works well, but using the true adjoint appears to recover an image more similar to the original, unblurred image.


\section{Example 2: Blind Channel Estimation}

Another application that depends on a fast adjoint in the gradient computation is in some formulations of blind channel estimation.  In blind channel estimation, a single source sends an unknown signal over multiple channels with unknown responses.  Observers collect the output of each channel and collectively attempt to determine the source signal and the channel impulse response from each channel.  Let $s$ be the unknown source signal of length $N$ and $h_i$ the channel impulse response of the $i^\text{th}$ channel, each of length $K$.  Then the output of the $i^\text{th}$ channel is

\[ x_i = h_i\ast s, \] 

\noindent where $\ast$ denotes a linear convolution.  Here we prescribe zero-padded boundary conditions, and take $x_i$ to be the ``full'' convolution of length $K+N-1$.  We hope to recover the source $s$ and channel responses $h_i$ from the output signals $x_i$.   However, note that there is both a magnitude and phase ambiguity in $h_i$ and $s$, since we can multiply $h_i$ by $\alpha\neq 0$ and divide $s$ by $\alpha$, leaving $x_i$ unchanged.

For simplicity, assume we have a single channel $h$ with observed output $x$; we will extend the problem to more channels in an example to follow.  We can pose the blind deconvolution problem for a single channel as

\begin{equation}
\label{eq:simple_bce_problem}
\min_{h,s}~ \frac{1}{2}\|h\ast s-x\|_2^2 + \lambda_h\|h\|_1 + \lambda_s \|s\|_1,
\end{equation}

\noindent where we include the terms $\lambda_h\|h\|_1$ and $\lambda_s\|s\|_1$ to promote sparsity in the recovered signals.  To promote other structure in the recovered signals (assuming it is present in the unknown true signals), we can include other terms, such as $\|h\|_\text{TV}$, the 1d total-variation (TV) seminorm.

Since (\ref{eq:simple_bce_problem}) involves $h\ast s$, it involves products of the problem variables and so it is a non-convex problem.  However, note that the outer product matrix $hs^T$ contains all of the terms used in the computation of the convolution $h\ast s$.  Indeed, we may define a \emph{linear} operator on $\mathcal{A}(hs^T)$ that computes the convolution: $\mathcal{A}(hs^T)=h\ast s$ \cite{ahmed_2013}.  Of course this is still a nonlinear operation, but the nonlinear operation is confined to the outer product $hs^T$.  In the following gradient computations, we will implicitly be computing the adjoint $\mathcal{A}^\ast$.

Let $f(h,s)=\frac{1}{2}\|h\ast s -x\|_2^2$ be the differentiable part of the objective.  First order methods require $\nabla_hf$ and $\nabla_sf$.  Note that when computing $\nabla_hf$, we hold $s$ fixed, and in doing so we may view the convolution $h\ast s$ as a linear operation on $h$.  Similarly, when computing $\nabla_s f$, we hold $h$ fixed.  Using the results of \cite{claerbout_1992}, we may determine a fast method for evaluating these gradients.  To illustrate the derivation, consider the case $K=3$ and $N=5$.  Set $x_\text{est}=h\ast s$.  We can write the convolution as a matrix-vector product in a couple different ways.  The first,

\[ x_\text{est} = \begin{bmatrix} h[0]\\h[1]&h[0]\\h[2]&h[1]&h[0]\\&h[2]&h[1]&h[0]\\&&h[2]&h[1]&h[0]\\&&&h[2]&h[1]\\&&&&h[2]\end{bmatrix}\begin{bmatrix}s[0]\\s[1]\\s[2]\\s[3]\\s[4]\end{bmatrix}, \] 

\noindent forms the matrix with entries from $h$, and the second,

\[ x_\text{est} = \begin{bmatrix} s[0]\\s[1]&s[0]\\s[2]&s[1]&s[0]\\s[3]&s[2]&s[1]\\s[4]&s[3]&s[2]\\&s[4]&s[3]\\&&s[4]\end{bmatrix}\begin{bmatrix}h[0]\\h[1]\\h[2]\end{bmatrix}
= S\, h, \] 

\noindent uses the entries of $s$ in the matrix $S$.  Note that here $h\ast s$ computes the ``full'' convolution, which includes using the zero-padded boundaries.  In {\sc matlab} we can compute the convolution with the \verb|conv| function:

\begin{verb}
  x_est = conv(h, s, 'full')
\end{verb}

For the gradient $\nabla_hf$, we treat $s$ as a constant, and so it is natural to use the second matrix-vector product where we use the entries of $s$ to form the linear operator.  We have

\[ \nabla_hf = \nabla_h\frac{1}{2}\|Sh-x\|^2 = S^\ast(Sh-x) = S^\ast r, \] 

\noindent where we define the residual $r=Sh-x = x_\text{est}-x$.  Written out, this gradient is

\[ \nabla_h f = \begin{bmatrix} \bar{s}[0]&\bar{s}[1]&\bar{s}[2]&\bar{s}[3]&\bar{s}[4]\\&\bar{s}[0]&\bar{s}[1]&\bar{s}[2]&\bar{s}[3]&\bar{s}[4]\\&&\bar{s}[0]&\bar{s}[1]&\bar{s}[2]&\bar{s}[3]&\bar{s}[4]\\\end{bmatrix}\begin{bmatrix}r[0]\\r[1]\\r[2]\\r[3]\\r[4]\\r[5]\\r[6]\end{bmatrix}. \] 

\noindent We recognize this as the cross-correlation between $s$ and $r$:

\[ \nabla_h f[n] = \sum_{k=0}^4 \bar{s}[k]r[k+n] = \big(\bar{s}[{-\,\cdot}]\ast r[\,\cdot\,]\big)[n], \]

\noindent where $\bar{s}[{-\,\cdot}]$ is the matched filter of $s[\,\cdot\,]$.  Note that the convolution here does not include the entries that use the zero-padded boundaries.  In {\sc matlab}, we can compute this convolution with

\begin{verb}
  Dfh = conv(r, conj(s(end:-1:1)), 'valid')
\end{verb}

\noindent The derivation and computations for $\nabla_sf$ are very similar.

Now that we can efficiently compute the required gradients, we can use existing first-order methods.  As an example, we consider a simulated underwater acoustic channel.  A single unknown real source signal is broadcast over two noisy acoustic channels with unknown real impulse responses.  We extend the blind channel estimation problem to two channels quite naturally:

\begin{align*}
   \min_{h_i,s} ~&\sum_{i=1}^2\left(\left\|h_i\ast s - x_i\right\|_2^2 + \lambda_{h}\|h_i\|_1\right) + \lambda_s\|s\|_1
\end{align*}

\noindent Further, we may add a differentiable 1d total-variation-like term in order to promote sharp transitions and stretches of nearly constant signal values (e.g., where the channel impulse response is constant for brief periods).  The true total-variation seminorm of $h$ may be computed as $\|\mathcal{D}h\|_1$, where $\mathcal{D}$ is the operator that subtracts consecutive pairs of elements of $h$: $\{\mathcal{D}h\}_j = h_{j+1} - h_j$.  We can in fact form $\mathcal{D}$ as a sparse matrix and compute $\mathcal{D}^\ast$ easily.  To ``soften'' the TV norm to a differentiable term, we may use the Huber function, defined as

\[ L_\delta(z) = \left\{\begin{array}{ll} \frac{1}{2}z^2 & |z| \le \delta\\ \delta\left(|z| - \frac{1}{2}\delta\right) & |z| > \delta. \end{array}\right. \] 

\noindent The Huber function is commonly used in regression with outliers, as it is linear for $|z|>\delta$, and so it is less sensitive to outliers.  Unlike ${\|\cdot\|_1}$, the Huber function is differentiable for all $z$.  We can approximate the TV norm with ${\|h\|_\text{TV} \defeq\|\mathcal{D}h\|_1 \approx L_\delta(\mathcal{D}h)/\delta}$.  The augmented problem, with the differentiable total-variation-like term is

\begin{align}
   \label{eq:bce_aug}
   f_\delta(h_i,s) &= \|h_i\ast s - x_i\|_2^2 + \lambda_{h,\text{TV}}L_\delta(\mathcal{D}h_i)/\delta,\nonumber\\
   \min_{h_i,s} ~&\sum_{i=1}^2\left(f_\delta(h_i,s) + \lambda_h\|h_i\|_1\right) + \lambda_s\|s\|_1.
\end{align}

It is straightforward to include the ``soft'' TV term in the gradients, and we can use standard first-order algorithms to solve the augmented problem.  Simulated data containing an impulsive source, two underwater acoustic channel impulse responses, and the two outputs of the underwater channels were provided by \cite{rideout_2016}.  The simulation also added zero mean Gaussian noise with standard deviation $0.005$.  We used the two channel outputs to attempt to recover the impulsive source and channel impulse responses.

We take $\lambda_h = 0.1$, $\lambda_s=0.01$, $\lambda_{h,\text{TV}}=0.01$, and $\delta=0.1$ for the Huber functions $L_\delta(\cdot)$.  We initialized $h_i$ and $s$ with zero mean Gaussian random numbers with unit standard deviation.  In practice, we may not know the true lengths of $h_i$ and $s$, $K$ and $N$, respectively.  We do know the length of the observed signals, which is $K+N-1$.  In this example, $K=894$ and $N=1717$, and we guess $K_\text{est}=844$ and $N_\text{est}=1767$.  We used L1General \cite{schmidt_2010} to find a local solution of the augmented problem (\ref{eq:bce_aug}).  The initial guesses for $h_i$ and $s$ were drawn from the standard normal distribution.

\begin{figure}
   \centering
   \includegraphics[width=0.5\textwidth]{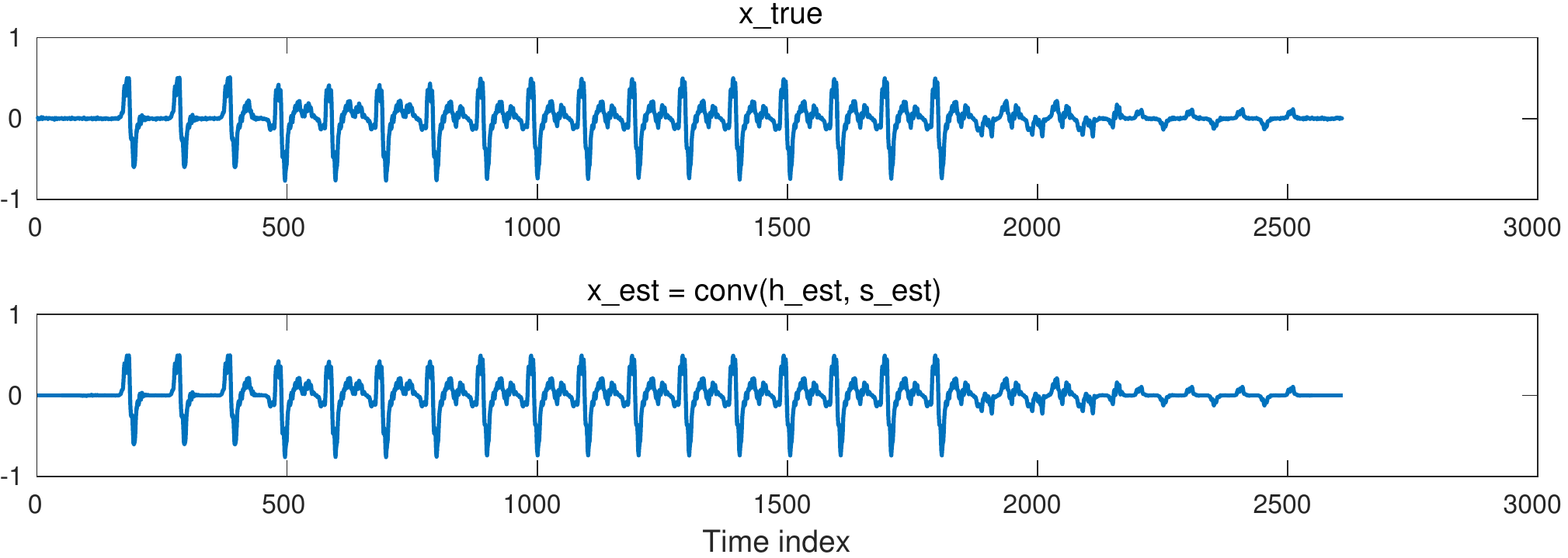} 
   \includegraphics[width=0.5\textwidth]{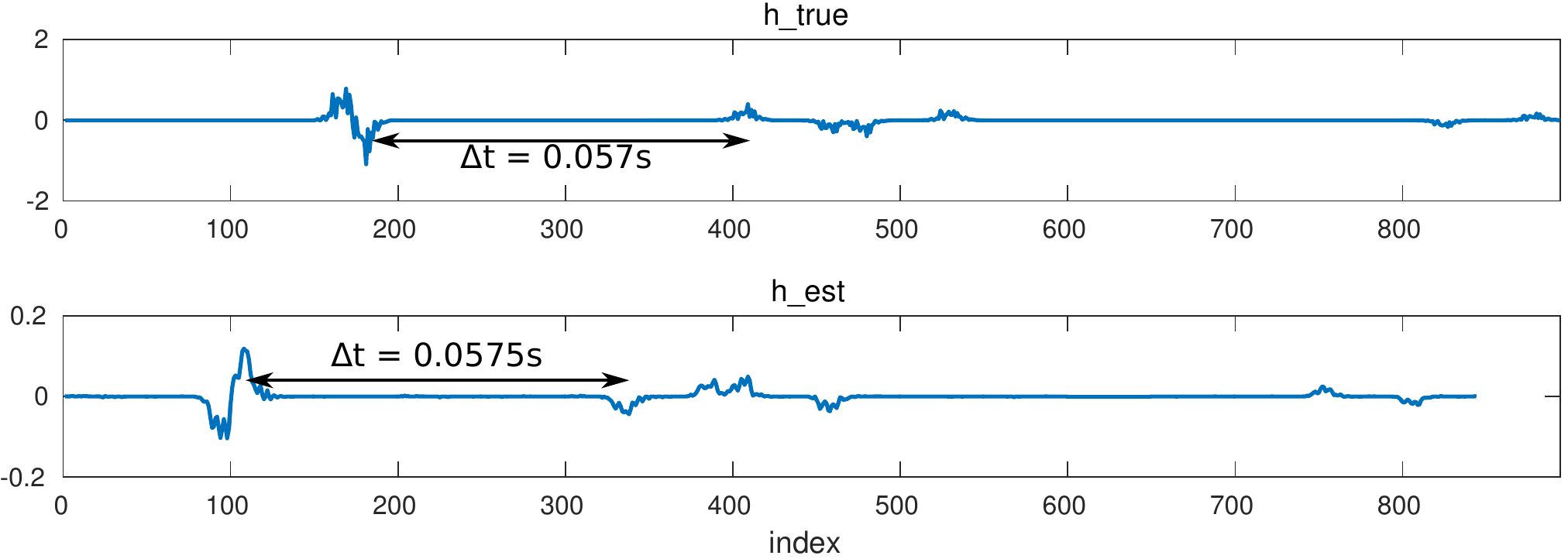}
   \includegraphics[width=0.5\textwidth]{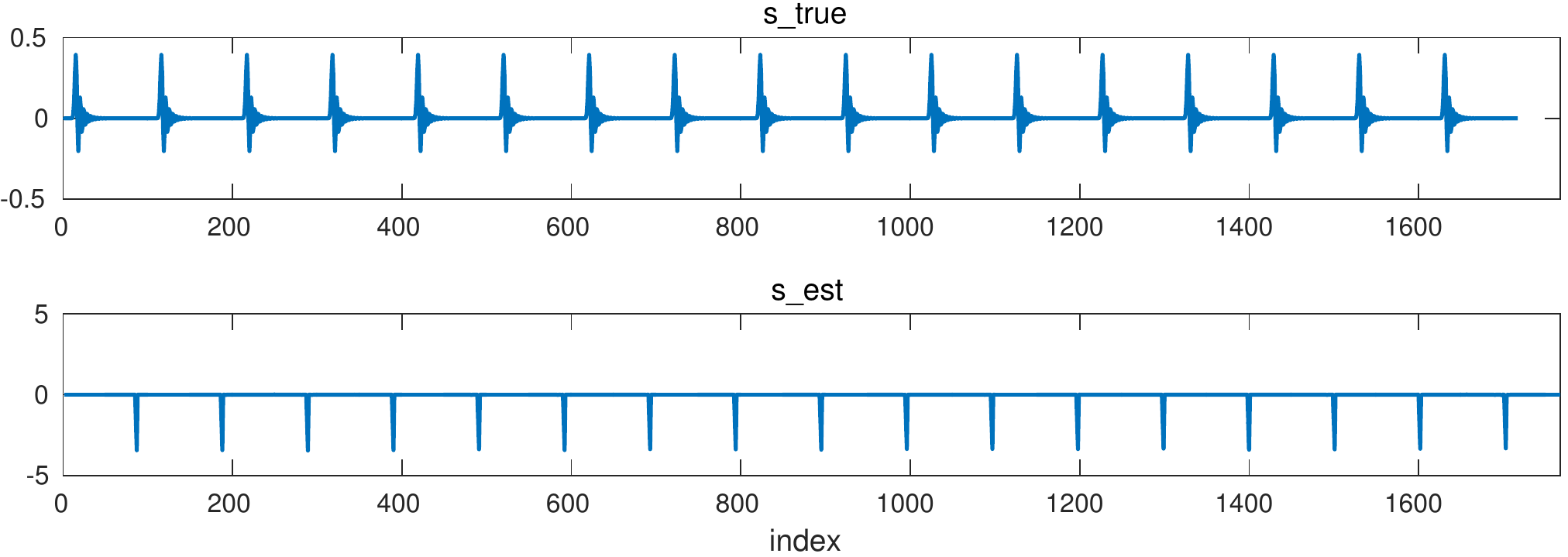}
   \caption{Blind channel estimation via problem (\ref{eq:bce_aug}) for a simulated underwater acoustic impulsive source.  Shown are the true and estimated first channel output, first channel impulse response, and source.  There are obvious magnitude and time shift differences, but the estimated channel output and channel impulse response are remarkably similar in structure and potentially useful for studying the underwater channel.}
   \label{fig:bce_rec}
\end{figure}

The results of the recovery are shown in Figure \ref{fig:bce_rec}.  We show the output and channel impulse response for only the first channel; the second channel is similar.  Researchers studying the acoustic channel are mainly interested in the time delays between peaks and relative phase shifts \cite{rideout_2016}.  Besides the overall magnitude ambiguity and time shift, it appears we can accurately estimate the time delays and phase shifts between the peaks of the channel response.

\section{Discussion}
We have demonstrated a framework for computing the adjoint of discrete wavelet transforms, which arise in gradient computations in optimization problems.  The framework separates the action of the wavelet transform and the signal extension operator.  The half-point symmetric extension is a common choice and we have provided its transpose, but other extension operators can readily be cast into the framework.  For many practical purposes, it appears that the adjoint operator is remarkably close to the pseudoinverse transform, however, implementing the true adjoint is just a simple matter of handling the signal extension operation.

Another interesting adjoint appears in a blind channel estimation optimization problem.  The convolution of two variables, which is a nonlinear operation, can be written as a matrix-vector product in two ways.  With this viewpoint, deriving the gradient becomes simple and a fast implementation is immediate.  For an impulsive source sent through a simulated underwater acoustic channel, the blind channel estimation problem is able to accurately recover important details of the acoustic channel response.


\bibliographystyle{IEEEtran}

\begin{thebibliography}{10}
\providecommand{\url}[1]{#1}
\csname url@samestyle\endcsname
\providecommand{\newblock}{\relax}
\providecommand{\bibinfo}[2]{#2}
\providecommand{\BIBentrySTDinterwordspacing}{\spaceskip=0pt\relax}
\providecommand{\BIBentryALTinterwordstretchfactor}{4}
\providecommand{\BIBentryALTinterwordspacing}{\spaceskip=\fontdimen2\font plus
\BIBentryALTinterwordstretchfactor\fontdimen3\font minus
  \fontdimen4\font\relax}
\providecommand{\BIBforeignlanguage}[2]{{%
\expandafter\ifx\csname l@#1\endcsname\relax
\typeout{** WARNING: IEEEtran.bst: No hyphenation pattern has been}%
\typeout{** loaded for the language `#1'. Using the pattern for}%
\typeout{** the default language instead.}%
\else
\language=\csname l@#1\endcsname
\fi
#2}}
\providecommand{\BIBdecl}{\relax}
\BIBdecl

\bibitem{hansen_2006}
P.~C. Hansen, J.~G. Nagy, and D.~P. O'Leary, \emph{Deblurring Images: matrices,
  spectra, and filtering}.\hskip 1em plus 0.5em minus 0.4em\relax SIAM, 2006.

\bibitem{beck_2009}
A.~Beck and M.~Teboulle, ``A fast iterative shrinkage-thresholding algorithm
  for linear inverse problems,'' \emph{SIAM Journal on Imaging Sciences}, 2009.

\bibitem{matlab_wt_2015}
\emph{MATLAB and Wavelet Toolbox Release R2015b}, The MathWorks, Inc., Natick,
  Massachusetts, 2015.

\bibitem{mallat_2009}
S.~Mallat, \emph{A Wavelet Tour of Signal Processing: The Sparse Way}.\hskip
  1em plus 0.5em minus 0.4em\relax Elsevier, 2009.

\bibitem{strang_1996}
G.~Strang and T.~Nguyen, \emph{Wavelets and Filter Banks}.\hskip 1em plus 0.5em
  minus 0.4em\relax Wellesley-Cambridge Press, 1996.

\bibitem{weber_1993}
G.~Weber, ``Usc-sipi image database: Version 4,'' 1993.

\bibitem{rout_2003}
S.~Rout, ``Orthogonal vs. biorthogonal wavelets for image compression,'' Ph.D.
  dissertation, Virginia Tech, 2003.

\bibitem{wang_2004}
Z.~Wang, A.~C. Bovik, H.~R. Sheikh, and E.~P. Simoncelli, ``Image quality
  assessment: from error visibility to structural similarity,'' \emph{Image
  Processing, IEEE Transactions on}, vol.~13, no.~4, pp. 600--612, 2004.

\bibitem{ahmed_2013}
A.~Ahmed, B.~Recht, and J.~Romberg, ``Blind deconvolution using convex
  programming,'' \emph{Information Theory, IEEE Transactions on}, vol.~60,
  no.~3, pp. 1711--1732, 2014.

\bibitem{claerbout_1992}
J.~F. Claerbout, \emph{Earth soundings analysis: Processing versus
  inversion}.\hskip 1em plus 0.5em minus 0.4em\relax Blackwell Scientific
  Publications Cambridge, Massachusetts, USA, 1992, vol.~6.

\bibitem{rideout_2016}
B.~Rideout and E.-M. Nosal, Personal communication, 2016.

\bibitem{schmidt_2010}
M.~Schmidt, ``Graphical model structure learning with l1-regularization,''
  Ph.D. dissertation, University of British Columbia (Vancouver), 2010.

\end{thebibliography}


\begin{IEEEbiographynophoto}{James Folberth}
(james.folberth@colorado.edu) is a Ph.D. student in Stephen Becker's information extraction group in the Department of Applied Mathematics at the University of Colorado at Boulder, USA.  His research interests include signal and image processing, machine learning, and optimization.
\end{IEEEbiographynophoto}

\begin{IEEEbiographynophoto}{Stephen Becker} is an assistant professor in the department of Applied Mathematics at  the University of Colorado at Boulder. He was the recipient of a postdoctoral fellowship to study at Paris 6, and the Goldstine postdoctoral fellowship at IBM Research. He completed his PhD in 2011 in Applied and Computational Mathematics from the California Institute of Technology, and his bachelors in Mathematics and Physics from Wesleyan University. His work on the TFOCS algorithm won the triennial Beale-Orchard-Hays Prize from the Mathematical Optimization Society.
   His current research is centered around deriving efficient optimization methods for large-scale data analysis, with emphasis on techniques for sparse and low-rank models.
\end{IEEEbiographynophoto}

This paper has supplementary downloadable material available at \url{https://github.com/jamesfolberth/ieee_adjoints}, provided by the authors. The material includes {\sc matlab} implementations of the adjoint operators discussed in this note. Contact \url{james.folberth@colorado.edu} for further questions about this work.

\vfill

\end{document}